\documentclass{article}
\usepackage{amsfonts}
\usepackage{graphicx}

\newcommand{\beq}{\begin{equation}}
\newcommand{\eeq}{\end{equation}}
\newcommand{\dpl}{\displaystyle}
\newcommand{\sech}{\mathop{\rm sech}\nolimits}

\begin{document}

\title{Yang-Baxter maps and integrable dynamics}

\maketitle

\begin{center}

{\bf A.P.Veselov }


\bigskip

{\it Department of Mathematical Sciences, Loughborough University,
Loughborough, Leicestershire, LE 11 3TU, UK
}

\bigskip

{\it Landau Institute for Theoretical Physics, Kosygina 2,

 Moscow, 117940,  Russia

\bigskip

e-mail: A.P.Veselov@lboro.ac.uk,
}

\end{center}

\bigskip

\bigskip

{\small  {\bf Abstract}
The hierarchy of commuting maps related to a set-theoretical solution of the quantum Yang-Baxter equation (Yang-Baxter map)
is introduced. They can be considered as dynamical analogues of the monodromy and/or transfer-matrices.
The general scheme of producing Yang-Baxter maps based on matrix factorisation is discussed in the context of the integrability problem
for the corresponding dynamical systems.
Some examples of birational Yang-Baxter maps coming from the theory of the periodic dressing chain and matrix KdV equation are discussed.}

\bigskip

\section*{Introduction}

In 1990 V.G. Drinfeld \cite{D} suggested to study the set-theoretical solutions
to the quantum Yang-Baxter equation. A. Weinstein and P. Xu \cite{WX} developed an approach to this problem 
based on the theory of Poisson Lie groups and symplectic groupoids.
P. Etingof, T. Schedler and A. Soloviev \cite{ESS}
(partly inspired by Hietarinta's work \cite{Hi}) made an extensive study of such solutions and 
the related algebraic and geometric structures (see also V.M. Buchstaber's paper \cite{Buch}).
Some interesting examples of such solutions have appeared in representation theory 
in relation with geometric crystals \cite{E}. 

While algebraic side of this problem has been fairly well understood
its dynamical aspects seem to be not appreciated yet. 

In this note we discuss the dynamical analogues of the monodromy and transfer-matrices which play
the crucial role in the theory of solvable models in statistical mechanics and quantum inverse scattering 
method
\cite{B, TF, JM}, where the quantum Yang-Baxter equation has its origin. 

The main property of the monodromy maps we introduce (which they share 
with the transfer-matrices) is the pairwise commutativity.
If the maps are polynomial or rational the commutativity relation is quite strong 
and sometimes implies the solvability of the corresponding dynamical system (see \cite{V} for the discussion of the known results in this direction). 
Unfortunately there are no general theorems known, so the question is whether this is true or not for 
the monodromy dynamics related to rational Yang-Baxter maps is open.

I am going to present some arguments in favour of the positive answer to this question. 
We will consider two classes of the set-theoretical solutions to the quantum Yang-Baxter equation coming from 
the theory of integrable systems.

The first one is related to matrix factorisations and QR-type of procedure. The simplest example here is the Adler's map describing the symmetries 
of the periodic dressing chain considered by A.B. Shabat and the author \cite{VS, A}. 

The second class is given by the interaction of solitons with the non-trivial internal parameters (e.g. matrix solitons). 
It is a well-known phenomenon in soliton theory (see e.g. \cite{Sol}, \cite{Kul}) that the interaction of $n$ solitons is completely determined by their pairwise 
interactions (so there are no multiparticle effects). The fact that the final result is independent of the order
of collisions means that the map determining the interaction of two solitons satisfies the Yang-Baxter relation.
We present an explicit form of this map for the matrix KdV equation using the formulas from \cite{G}.

\section*{Monodromy maps for the set-theoretical solutions of the quantum Yang-Baxter equation}

The original {\it quantum Yang-Baxter equation} is the following relation on a linear operator $R : V \otimes V \rightarrow V \otimes V:$
\begin{equation}
\label{YB}
R_{12} R_{13} R_{23} = R_{23} R_{13} R_{12},
\end{equation}
where $R_{ij}$ is acting in $i$-th and $j$-th components of the tensor product $V \otimes V \otimes V$
(see e.g. \cite{JM}).

Following Drinfeld's suggestion \cite{D} let us consider the following set-theoretical version of this equation.

Let $X$ be any set, $R: X \times X \rightarrow X \times X$ be a map from its square into itself.
Let $R_{ij}: X^{n} \rightarrow X^{n}, \quad X^{n} = X \times X \times .....\times X$
be the maps which acts as $R$ on $i$-th and $j$-th factors and identically on the others.
More precisely, if $R(x, y) = (f(x,y), g(x,y)), x,y \in X$ then 
$R_{ij} (x_1, x_2, \dots, x_n) = (x_1, x_2, \dots, x_{i-1}, f(x_i,x_j), x_{i+1}, \dots, x_{j-1}, g(x_i,x_j), x_{j+1},
\dots,x_n)$
for $i<j$ and 
$(x_1, x_2, \dots, x_{j-1}, g(x_i,x_j), x_{j+1},\dots, x_{i-1}, f(x_i,x_j), x_{i+1},\dots,x_n)$
otherwise. 
In particular for $n=2$
$R_{21}(x,y) = (g(y,x), f(y,x)).$ If $P: X^2 \rightarrow X^2$ is the permutation of $x$ and $y$: $P(x,y) = (y,x)$, then
obviously we have $$R_{21} = P R P.$$ 

{\bf Definition.} We will call $R$ {\it Yang-Baxter map} if it satisfies the Yang-Baxter relation (\ref{YB}) 
considered as the equality of the maps of $X \times X \times X$ into itself. 
If additionally $R$ satisfies the relation
\begin{equation}
\label{U}
R_{21} R = Id,
\end{equation}
we will call it {\it reversible Yang-Baxter map}.

The condition (\ref{U}) is usually called the {\it unitarity} condition, but in our case
the term {\it reversibility} is more appropriate since in dynamical systems terminology 
this condition means that the map $R$ 
is reversible with respect to the permutation $P$. 

{\bf Remark.} If we consider the linear space $V = {\bf C}^{X}$ spanned by the set $X$, then any Yang-Baxter map $R$
induces a linear operator in $V \otimes V$ which satisfies the quantum Yang-Baxter equation in the usual sense.
Therefore we indeed have a (very special) class of solutions to this equation. 
If $X$ is a finite set and $R$ is a bijection then we have permutation-type solutions discussed in \cite{Hi}.
However this point of view on the Yang-Baxter maps seems to be artificial because it does not reflect the nature of the maps 
(like birationality). That's why I prefer the term "Yang-Baxter map"
rather than "set-theoretical solution to the quantum Yang-Baxter equation" (which also has a disadvantage of being too long).
\footnote {P. Etingof \cite{E} suggested (in birational situation) the alternative term "rational set-theoretical $R$-matrix" 
which also is not ideal: it does not make clear that we are talking about maps.}

We can represent the relations (\ref{YB}, \ref{U}) in the standard diagrammatic way as follows:

\begin{center}

\includegraphics[width=10cm]{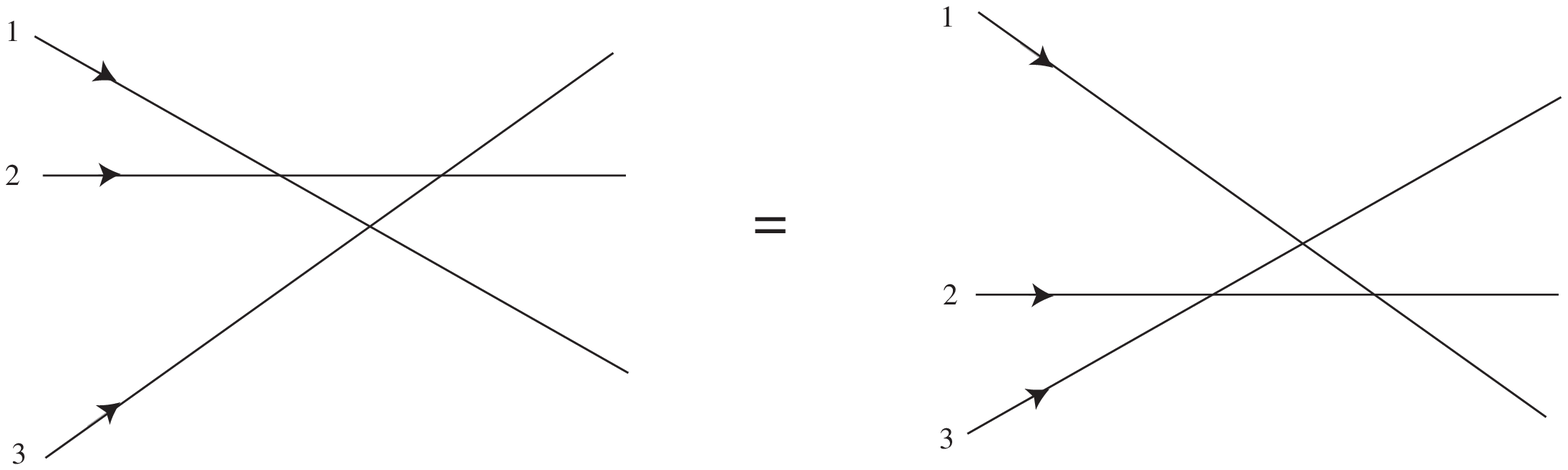}

Figure 1

\includegraphics[width=10cm]{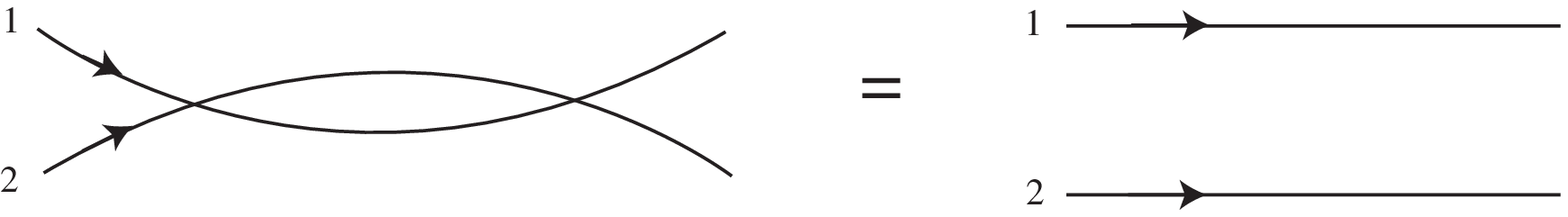}

Figure 2

\end{center}

One can introduce also a more general parameter-dependent Yang-Baxter
equation as the relation
\begin{equation}
\label{sYB}
R_{12}(\lambda_1, \lambda_2) R_{13}(\lambda_1, \lambda_3) R_{23} (\lambda_2, \lambda_3) = R_{23}(\lambda_2, \lambda_3)
R_{13}(\lambda_1, \lambda_3) R_{12}(\lambda_1, \lambda_2)
\end{equation}
and the corresponding unitarity (reversibility) condition as
\begin{equation}
\label{sU}
R_{21}(\mu,\lambda) R(\lambda, \mu) = Id.
\end{equation}
Although it can be reduced to the usual case
by considering $\tilde X = X \times {\bf C}$ and $\tilde R (x,\lambda; y, \mu) = R (\lambda,\mu) (x,y)$
sometimes it is useful to keep the parameter separately (see the examples below).

In the theory of the quantum Yang-Baxter equation the following {\it transfer-matrices} 
$t^{(n)}: V^{\otimes n} \rightarrow V^{\otimes n}$ play a fundamental role.
They are defined as the trace of the monodromy matrix 
$$t^{(n)} = tr_{V_0} R_{0 n} R_{0 n-1} \dots R_{01}$$ with respect to the additional space $V_0$.
If the solution of the quantum Yang-Baxter equation depends on an additional spectral parameter 
$\lambda$ then the corresponding transfer-matrices commute:
$$t^{(n)}(\lambda) t^{(n)}(\mu) = t^{(n)}(\mu) t^{(n)}(\lambda),$$
which is the crucial fact in that theory (see \cite{B,TF,JM}).

In the general set-theoretical case we have a problem with the trace operation, so it seems that we do not have
a direct analogue of the transfer-matrices.
Our main idea is to replace them by the following maps, which combine the properties of both monodromy and transfer-matrices.
\footnote{As P.P. Kulish explained to me at NEEDS conference 
(Cadiz, June 2002) the corresponding operators (sometimes called as Yang's operators) are well-known to the experts in the theory of Bethe ansatz.
They first appeared in the C.N.Yang's paper \cite{Yang} and play an essential role in Frenkel - Reshetikhin construction
of the q-KZ equation \cite{FR} and in Fomin - Kirillov approach to the theory of Schubert polynomials and related symmetric functions \cite{FK}.
The best discussion of these operators (with explanation of their relations to transfer-matrices) I have found in M. Gaudin's book \cite{Gaudin} 
(see especially Chapter 10, sections 2 and 3).}

Let us define the {\it monodromy maps} $T_i^{(n)}, i = 1,\dots, n$ as the maps of $X^{n}$ into 
itself by the following formulas:
\begin{equation}
\label{T}
T_i^{(n)} = R_{i i+n-1} R_{i i+n-2} \dots R_{i i+1},
\end{equation}
where the indeces are considered modulo $n$ with the agreement that we are using $n$ rather than 0.
In particular $T_1^{(n)} = R_{1 n} R_{1 n-1} \dots R_{12}.$

{\bf Theorem 1.} {\it For any reversible Yang-Baxter map} $R$ {\it the monodromy maps} $T_i^{(n)},\quad i =1,\dots, n$ {\it commute with each other:}
\begin{equation}
\label{comm}
T_i^{(n)} T_j^{(n)} = T_j^{(n)} T_i^{(n)}
\end{equation}
{\it and satisfy the property} 
\begin{equation}
\label{prod}
T_1^{(n)} T_2^{(n)} \dots T_n^{(n)} = Id.
\end{equation}
{\it Conversly, suppose that the maps} $T_i^{(n)}$ {\it determined by the formula (\ref{T}) commute
and satisfy the relation (\ref{prod}) for any $n \geq 2$ then $R$ is a reversible Yang-Baxter map.}

Proof of the first part is very similar to the proof of the commutativity of the transfer-matrices 
and follows from the consideration of the corresponding diagrams
representing the products $T_i^{(n)} T_j^{(n)}$ and $T_j^{(n)} T_i^{(n)}$ respectively (cf \cite{JM}):
\begin{center}
\includegraphics[width=10cm]{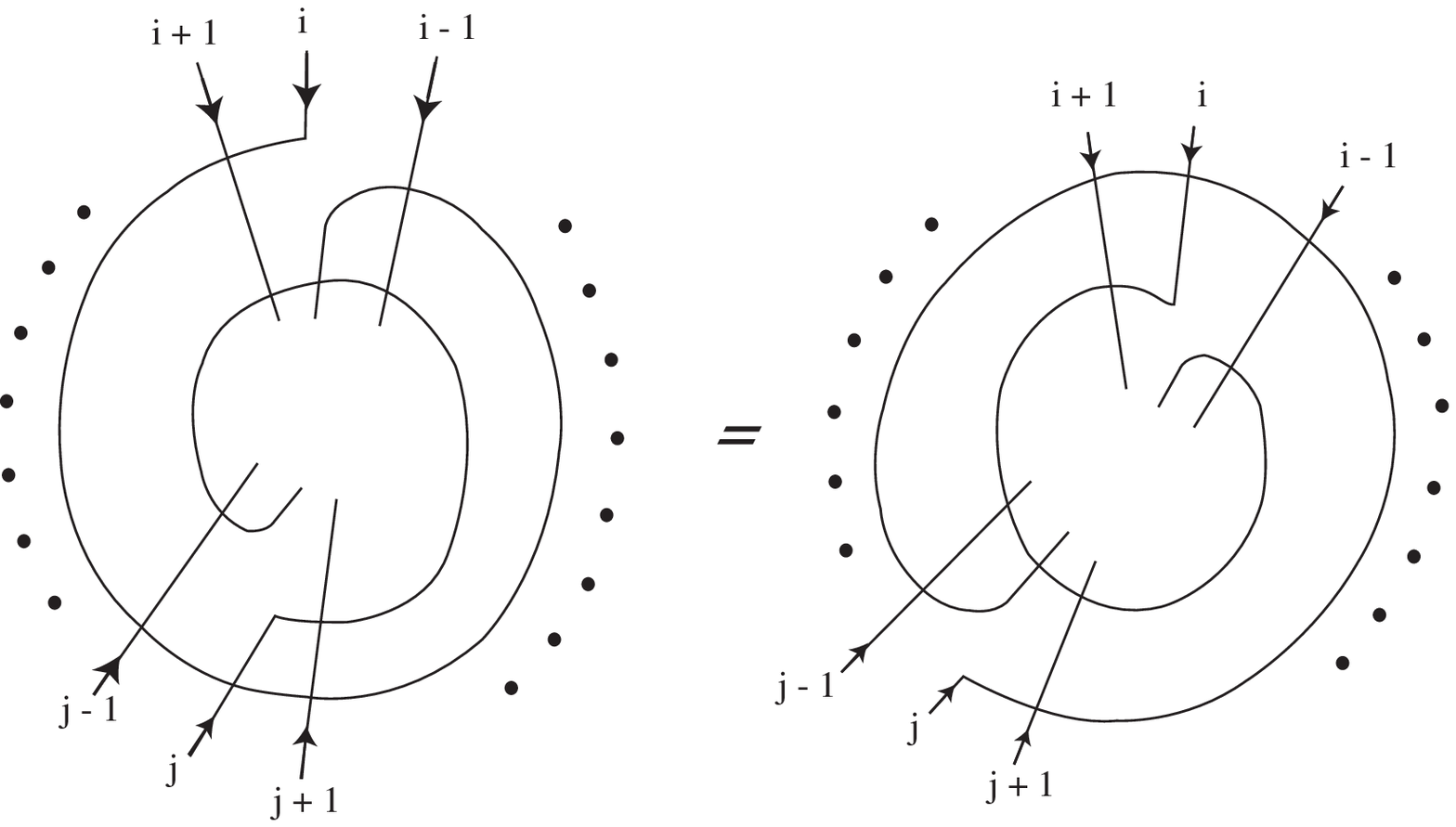}

Figure 3
\end{center}


One can easily check that the second diagram is a result of the several operations presented
above on figures 1 and 2 applied to the first diagram.
The second identity also follows from similar consideration (see fig. 4 below where the case $n=3$ is presented).

\begin{center}
\includegraphics[width=10cm]{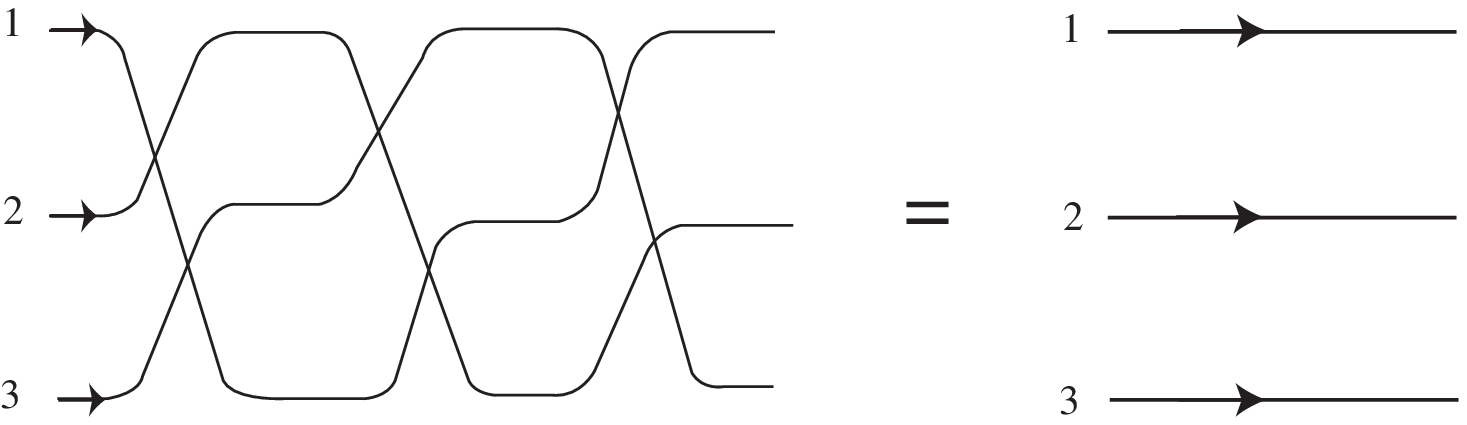}

Figure 4
\end{center}

To prove the second part it is actually enough to consider only the cases $n=2$ and $n=3.$
When $n=2$ we have two maps $T_1^{(2)} = R$ and $T_2^{(2)} = R_{21}$, so the relation (\ref{prod})
becomes simply the unitarity (reversibility) condition (\ref{U}). The commutativity condition is trivial in this case.
For $n=3$ then we have three transfer-maps: $T_1^{(3)} = R_{13} R_{12}, T_2^{(3)} = R_{21} R_{23}, T_3^{(3)} = R_{32} R_{31}.$
The product $T_1^{(3)} T_2^{(3)} T_3^{(3)} =  R_{13} R_{12} R_{21} R_{23} R_{32} R_{31}$ is obviously an identity because $R_{12} R_{21}=R_{23} R_{32}=R_{13} R_{31} =Id.$
Commutativity $T_1^{(3)} T_2^{(3)} = T_2^{(3)} T_1^{(3)}$ means that $R_{13} R_{12} R_{21} R_{23} = R_{21} R_{23} R_{13} R_{12}$ which is equivalent to
Yang-Baxter equation (\ref{YB}) modulo the unitarity relation which was already shown. This proves the theorem.

The properties of the monodromy maps can also be explained in terms of the corresponding extended affine Weyl group $\tilde A_{n-1}^{(1)}.$
Indeed it is well-known that the maps $$S_{i} = P_{ii+1} R_{ii+1},$$
where $P_{ij}$ is the permutation of $i$-th and $j$-th factors in $X^{n},$
satisfy the relations
\begin{equation}
\label{aff}
S_{i}S_{i+1}S_{i} = S_{i+1}S_{i}S_{i+1}, \quad  i=1, \dots, n
\end{equation}
and
\begin{equation}
\label{inv}
S_i^2 = Id
\end{equation}
which are the defining relations of the affine Weyl group $A_{n-1}^{(1)}.$
One can check that the monodromy maps $T_i^{(n)}$ correspond to the translations
in the extended affine Weyl group $\tilde A_{n-1}^{(1)}$ (which is the normalisator of $A_{n-1}^{(1)}$) 
generated by $S_1, \dots, S_n$ and
the cyclic permutation $\omega = P_{1 n} P_{1 n-1} \dots P_{12},\quad \omega^{n} = Id$
(see \cite{Bour}).

{\bf Remark.} If we consider only the set of $S_i$ with $i= 1, \dots, n-1$ (without $S_n$) 
then they generate the (finite) permutation group ${\bf S}_{n}.$ In that case as it was shown by Etingof et al
\cite{ESS} the action corresponding to any (non-degenerate) Yang-Baxter map is equivalent to the usual action
of ${\bf S}_{n}$ by permutations (and thus trivial from dynamical point of view).
It is very important that in our case the monodromy maps generate the actions of an infinite
abelian group ${\bf Z}^{n-1}$ so the question of whether two such actions are equivalent or not 
is a non-trivial problem of dynamical systems theory.

In the next two sections we present some examples of the Yang-Baxter maps and 
explain how they appear in the theory of integrable systems.
In all the examples we present the Yang-Baxter maps are birational and depend on the spectral parameters,
so we consider the relations (\ref{sYB}, \ref{sU}) rather then (\ref{YB}, \ref{U}).

\section*{Yang-Baxter maps and matrix factorisations}

Let $A(x,\zeta)$ be a family of matrices depending on the point $x \in X$ and  
a "spectral" parameter $\zeta \in {\bf C}$. One should think of $X$ being an algebraic variety
and $A$ depending polynomially/rationally on $\zeta$ although the procedure works always when
the factorisation problems appeared are uniquely solvable. For matrix polynomials usually this is the case
once the factorisation of the corresponding determinant is fixed (see \cite{GLR}). 

The following procedure is a version of standard QR-algorithm in linear algebra and known to be very useful 
in discrete integrable systems \cite{MV}. 

Consider the product $L = A(x,\zeta) A(y,\zeta),$ then change the order
of the factors $L \rightarrow \tilde{L} = A(y, \zeta) A(x, \zeta)$ and refactorise it again:
$\tilde{L} = A(\tilde{x}, \zeta) A(\tilde{y}, \zeta).$
Now the map $R$ is defined by the formula 
\beq
\label{fact}
R(x,y) = (\tilde{x},\tilde{y}).
\eeq

The first claim is that this map satisfies the relations (\ref{YB}, \ref{U}). \footnote{I was not able to find out to whom
this simple but important observation should be prescribed to. It probably goes back to A.B. Zamolodchikov's works in the late 70-th on 
factorised S-matrices and related algebras. I would like to mention in this relation also the papers \cite{RS}, \cite{BFZ} and the recent preprint \cite{O}.
The first concrete example I know which came from this construction is the Adler's map considered below.}
Indeed if we consider the product $A(x_3) A(x_2) A(x_1)$ (we omit here the $\zeta$ for shortness) then applying
the left hand side of (\ref{YB}) to this product we have
$A(x_3) A(x_2) A(x_1) = A(x_3^{(1)}) A(x_1^{(1)}) A(x_2^{(1)})=
A(x_1^{(2)}) A(x_3^{(2)}) A(x_2^{(2)}) = A(x_1^{(3)}) A(x_2^{(3)}) A(x_3^{(3)}).$
Similarly the right hand side corresponds to the relations $A(x_3) A(x_2) A(x_1) = \\
A(\tilde{x}_2^{(1)}) A(\tilde{x}_3^{(1)}) A(\tilde{x}_1^{(1)})=
A(\tilde{x}_2^{(2)}) A(\tilde{x}_1^{(2)}) A(\tilde{x}_3^{(2)}) = 
A(\tilde{x}_1^{(3)}) A(\tilde{x}_2^{(3)}) A(\tilde{x}_3^{(3)}).$
Because of the uniqueness of the factorisation we have $x_i^{(3)} = \tilde{x}_i^{(3)},$ which is the Yang-Baxter relation.

The unitarity relation (\ref{U}) is obvious.

Let us consider now the corresponding monodromy maps $T_i^{(n)}$. We claim that they have many integrals.
Indeed let us introduce the {\it monodromy matrix} $$M= A(x_n,\zeta) A(x_{n-1},\zeta) \ldots A(x_1, \zeta).$$

{\bf Theorem 2.} {\it The monodromy maps} $T_i^{(n)}, \quad i=1,\ldots, n$ {\it related to the Yang-Baxter map} (\ref{fact}) 
{\it preserve the spectrum of the corresponding monodromy matrix} $M(x_1, \ldots, x_n; \zeta)$ {\it for all} $\zeta.$

The proof is simple. For the map $T_1^{(n)} = R_{1 n} R_{1 n-1} \dots R_{12}$ one can easily see that
$A(x_n,\zeta) A(x_{n-1},\zeta) \ldots A(x_1, \zeta) = A(\tilde{x}_1, \zeta)  A(\tilde{x}_n,\zeta) A(\tilde{x}_{n-1},\zeta)\ldots A(\tilde{x}_2, \zeta),$
where $(\tilde{x}_1, \ldots ,\tilde{x}_n) = T_1^{(n)} (x_1, \ldots, x_n).$
Now the claim follows from the well-known fact that the spectrum of the product of the matrices is invariant under the cyclic permutation of the factors, so 
$$Spec M(x_1, \ldots, x_n; \zeta) = Spec M(\tilde{x}_1, \ldots ,\tilde{x}_n; \zeta).$$
To prove the same for $T_i^{(n)}$ one should consider the matrix $A(x_{n+i-1},\zeta) \ldots A(x_i, \zeta)$
which obviously has the same spectrum as $M$ and use the same arguments.

{\bf Corollary.} {\it The coefficients of the characteristic polynomial} $\xi = \det (M - \lambda I)$
{\it are the integrals of the commuting monodromy maps} $T_i^{(n)}.$

If the dependence of $A$ (and therefore $M$) on $\lambda$ is polynomial then the dynamics can be linearised
on the Jacobi varieties of the corresonding spectral curves $\xi(\zeta, \lambda) = 0$ (see \cite{MV}).
If these maps are also symplectic (which surprisingly happens quite often) then one can use also the discrete
version of the Liouville theorem \cite{V} to claim their integrability.

{\bf Remark.} Although the procedure is very similar to the one proposed in \cite{MV} for the discrete Lagrangian systems
there is one important difference: in the construction of Yang-Baxter maps the factors $A$ should depend only on one of the variables,
while in general scheme they usually depend both on $x$ and $y$ (see \cite{MV}).

{\bf Example}: Symmetries of the periodic dressing chain \cite{VS, A}.

Here $X = {\bf C}\times{\bf C}$ and the matrix $A(x), x = (f,\beta)$  has the form:

$$
\Bigl(
\begin{array}{cc}
f & 1\\
f^2 + \beta - \zeta & f
\end{array}
\Bigl)
$$

One can check that the factorisation procedure described above
leads in this case to the following birational map
 $R: (f_1,\beta_1; f_2, \beta_2) \rightarrow (\tilde{f_1},\beta_1; \tilde{f_2}, \beta_2):$ 
\begin{equation}
\label{Adler}
\tilde{f_1} = f_2 - \frac{\beta_1 - \beta_2}{f_1 + f_2} \qquad
\tilde{f_2} = f_1 - \frac{\beta_2 - \beta_1}{f_1 + f_2}.
\end{equation}
This map (modulo additional permutation) first appeared in Adler's paper \cite{A} as a symmetry of the periodic dressing chain
\cite {VS} and later was discussed in more details by Noumi and Yamada \cite{NY} (without reference to Adler's
paper \cite{A} which they seem to be not aware of). 

The fact that this map satisfies the Yang-Baxter equation and
unitarity condition (\ref{sYB}, \ref{sU}) follows from Adler's results (and of course, from the general claim above).
The monodromy maps in this example are symplectic and integrable in the sense of discrete Liouville 
theorem \cite{V} (see \cite{A}). 

V. Adler discovered also a remarkable geometric representation of these maps  
(recutting of polygons). It is interesting that the monodromy maps correspond to the cyclic recutting
procedure, which was of special interest for him (see \cite{A}).

\section*{Interaction of solitons as a Yang-Baxter map.}

There is another natural source of the Yang-Baxter maps in the theory of integrable systems:
soliton interaction. 

Consider any integrable by the inverse scattering method 
PDE in 1+1 dimensions which has multisoliton solutions. Suppose that each soliton has a non-trivial 
internal degrees of freedom described by the set $X$ (which is usually a manifold).
Then the standard arguments in soliton theory based on the existense of the commuting flows shows 
that the result of interaction of three solitons
must be independent of the order of collisions (see \cite{Kul}, \cite{G}) and thus the corresponding map $R$
is a set-theoretical solution to the quantum Yang-Baxter equation.

A good example is the matrix KdV equation:

$$
U_t+3UU_x+3U_xU+U_{xxx}=0,
$$
$U$ is $d \times d$ matrix.

It is easy to check that it has the soliton solution of the form $$U = 2 \lambda^2 P \sech^2 (\lambda x- 4\lambda^3 t),$$
where $P$ must be a projector: $P^2 = P.$ If we assume that $P$ has rank 1 than $P$ should have the form
$P=\frac{\dpl \xi\otimes \eta}{\dpl (\xi,\eta)}$. Here $\xi$ is a vector in a vector space $V$ of dimension $d$,
$\eta$ is a vector from the dual space $V^*$ (covector)
and bracket $(\xi,\eta)$ means the canonical pairing between $V$ and $V^*.$

The change of the matrix amplitudes $P$ ("polarisations") of two solitons with the velocities $\lambda_1$ and $\lambda_2$ 
after their interaction is described by the following map \cite{G}:
$$R(\lambda_1, \lambda_2): (\xi_1, \eta_1; \xi_2, \eta_2) \rightarrow (\tilde{\xi_1}, \tilde{\eta_1}; \tilde{\xi_2}, \tilde{\eta_2)}$$
\beq
\label{fi}
\tilde{\xi_1} = \xi_1+\frac{\dpl 2\lambda_2(\xi_1,\eta_2)}{\dpl (\lambda_1-
\lambda_2)(\xi_2,\eta_2)}\xi_2, \qquad
\tilde{\eta_1} = \eta_{1}+\frac{\dpl
2\lambda_2(\xi_2,\eta_1)}{\dpl (\lambda_1-
\lambda_2)(\xi_2,\eta_2)}\eta_2,
\eeq

\beq
\label{fj}
\tilde{\xi_2} = \xi_2+\frac{\dpl 2\lambda_1(\xi_2,\eta_1)}{\dpl (\lambda_2-
\lambda_1)(\xi_1,\eta_1)}\xi_1,\qquad
\tilde{\eta_2} = \eta_2+\frac{\dpl 2\lambda_1(\xi_1,\eta_2)}{\dpl (\lambda_2
-\lambda_1)(\xi_1,\eta_1)}\eta_1.
\eeq

In this example $X$ is the set of projectors $P$ of rank 1 which is the variety ${\bf CP}^{d-1} \times {\bf CP}^{d-1}.$
Although the integrability of the corresponding monodromy maps should be somehow a corollary of the
integrability of the initial PDE the precise meaning of this is still to be investigated 
(see \cite{GV} for the latest results in this direction).

I would like only to mention that the interaction map (\ref{fi}),(\ref{fj}) also can be described by the factorisation scheme if one takes
the matrices of the form $$A(P, \lambda; \zeta) = I + \frac{2 \lambda}{\zeta -\lambda} P$$
hinted by the inverse scattering problem for the matrix KdV equation.\footnote{Yuri Suris suggested a very simple explanation of 
this form of the matrix $A$ which works for Adler's maps as well. The question how far this observation goes is currently under investigation \cite{SV}.}

\section*{Discussion}

We have shown that with any Yang-Baxter map one can relate the hierarchy of commuting
monodromy maps $T_i^{(n)}, \quad n=2,3, \ldots.$ In the case when the maps are symplectic and have enough integrals in involution
one can use the discrete Liouville theorem \cite{V} to claim the integrability.
If these maps are related to matrix factorisation problems then one can use the Abel map to linearise the dynamics 
on the jacobians of the spectral curves \cite{MV}. 

But for the general maps there are no natural definitions of integrability.
For algebraic maps one can use the symmetry approach based on the existense of the commuting maps (see \cite{V}). 
In particular, it is known after G. Julia and P. Fatou that if two polynomial maps of $\bf{C}$ into itself $p: z \rightarrow p(z)$
and $q: z \rightarrow q(z)$ commute and have no common iterations then they must be equivalent either to the powers $z^k$ or
to Chebyshev maps $T_k(z)$ (see \cite{V} for the references and further developments in this direction).
In both cases the dynamics can be described by explicit formulas.

This means that one should expect some sort of integrability for the dynamics determined by the monodromy maps
(at least for the birational Yang-Baxter maps). To make this precise is the main challenge in this area.
Probably one should start with the investigation of the Yang-Baxter maps, 
which are polynomial automorphisms ( or, more general 
birational Cremona maps) of ${\bf C}^2 = {\bf C} \times {\bf C}$. 

I would like to mention in this relation the following observation due to V. Lyubashenko:
in case when map $R$ is of the form $R(x,y) = (p(x),q(y))$ the Yang-Baxter relation (\ref{YB}) is equivalent to the commutativity
of the maps $p$ and $q$ (the unitarity relation (\ref{U}) in this case simply means that $q = p^{-1}$).
It is interesting that this was one of the examples which stimulated V. Drinfeld to raise the general
question about the set-theoretical solutions to the quantum Yang-Baxter equation (see \cite{D}).

\section*{ Acknowledgements}
My initial interest to the Yang-Baxter maps has been originated
in the discussions with V.M. Goncharenko of the interaction of matrix solitons.
I am grateful to V.G. Drinfeld who introduced me into the history of the 
problem and to P. Etingof, S. Fomin, V. Papageorgiou, E. Sklyanin, F. Smirnov and M. Sevryuk for useful discussions
and comments. 

I am grateful also to the organisers and participants of the NEEDS conference in Cadiz (10-15 June 2002) 
and SIDE-V conference in Giens (21-26 June 2002) and especially to P.P. Kulish, A.B. Shabat and Yu. Suris 
for stimulating and helpful discussions.

I would like also to thank Helen Sherwood who helped me with preparation of the pictures.

\end{document}